\documentclass[a4paper, 12pt]{amsart}
\numberwithin{equation}{section}
\usepackage{amssymb}
\usepackage{amsmath}
\usepackage{amsfonts}
\usepackage{marvosym}
\usepackage{tikz}
\newcommand*\circled[1]{\tikz[baseline=(char.base)]
{\node[shape=circle,draw,inner sep=1] (char) {#1};}}
\usepackage{stackrel}

\newtheorem{theorem}{Theorem}[section]
\newtheorem{thm}{Theorem}[section]
\newtheorem{lem}[theorem]{Lemma}
\newtheorem{cor}[theorem]{Corollary}
\newtheorem{prop}[theorem]{Proposition}

\theoremstyle{definition}

\newtheorem{rem}[theorem]{Remark}

\def\cb{{\mathcal B}}

\def\ce{{\mathcal E}}

\def\ch{{\mathcal H}}

\def\cs{{\mathcal S}}

\def\ga{{\mathfrak A}}
\def\gb{{\mathfrak B}}

\def\bc{{\mathbb C}}

\def\bm{{\mathbb M}}
\def\bn{{\mathbb N}}
\def\bp{{\mathbb P}}

\def\bz{{\mathbb Z}}

\def\a{\alpha}
\def\b{\beta}
  
  \def\D{\Delta}
\def\eeps{\epsilon}
\def\eps{\varepsilon}
\def\io{\iota}

\def\p{\pi}
\def\n{\nu}

 \def\S{\Sigma}

\def\f{\varphi}  \def\F{\Phi}
\def\th{\theta} 
\def\om{\omega}

\def\id{\hbox{id}}

\def\aut{\mathop{\rm aut}}

\newcommand{\ty}[1]{\mathop{\rm {#1}}}
\def\di{{\rm d}}
\def\id{{\rm id}}

\def\ad{\mathop{\rm ad}}

\def\idd{{1}\!\!{\rm I}}

\newcommand{\nn}{\nonumber}

\DeclareMathAlphabet{\mathpzc}{OT1}{pzc}{m}{it}

\begin{document}

\title[de finetti theorem]
{Symmetric states for $C^*$-Fermi systems I: de finetti theorem}
\author{Francesco Fidaleo}
\address{Francesco Fidaleo\\
Dipartimento di Matematica \\
Universit\`{a} di Roma Tor Vergata\\
Via della Ricerca Scientifica 1, Roma 00133, Italy} \email{{\tt
fidaleo@mat.uniroma2.it}}
\date{\today}

\begin{abstract}
In the present note, which is the first part of a work concerning the study of the set of the symmetric states for Fermi systems,
we describe the extension of the De Finetti theorem to the infinite Fermi $C^*$-tensor product of a single (separable) general $\bz^2$-graded $C^*$-algebra.
\end{abstract}

\maketitle

\section{Introduction}
\label{sec1}

The investigation of the exchangeable stochastic processes, or the connected problem of the symmetric probability measures, is an important tool in probability theory. In this context, De Finetti Theorem deals with the relationship between independence and exchangeability.

More precisely, let $X_1, X_2,\dots$ be a countably infinite sequence of random variables. Such a sequence is said to be exchangeable if, for each integer $n$ and any choice of a finite subset $X_{i_1},\dots,X_{i_n}$ of $n$ distinct random variables, it has the same joint probability distribution of any other choice of $n$ distinct random variables $X_{j_1},\dots,X_{j_n}$. 

Obviously, if the sequence $X_1, X_2,\dots$ is independent and identically distributed, it is exchangeable, but the converse is false. De Finetti theorem simply asserts that an exchangeable sequence is indeed a "mixture" of independent and identically distributed sequences of random variables in a sense which is clarified below.\footnote{The original form of De Finetti theorem was concerning sequences of Bernoulli random variables, a natural generalisation to quite general sequences was provided by Hewitt and Savage in \cite{HS}.}

Due to the application to a an enormity of cases of interest, De Finetti theorem, firstly appeared in \cite{DeF}, represents a milestone in probability theory, and it has many interesting equivalent formulations. It should also be remarked that several attempts to generalise those De Finetti-like results to the quantum situation have been done in the last years. Since the aim of the present note is to provide the extension of the De Finetti theorem to quantum models describing  the presence of Fermi particles, among the early papers concerning some extension of such a fundamental theorem, we mention \cite{K}, and \cite{Fbo} for an understandable analysis on a toy-model.

The Hewitt-Savage formulation of the De Finetti theorem concerns the study of the class of probability Radon measures on the Tikhonov product
$\times_\bn X$ of a single compact metrisable space $X$ which are invariant under the action of the group of all permutations moving only finitely many indices. Such invariant probability measures are called {\it symmetric}. It can be immediately paraphrased in terms of the dual object consisting of the $C^*$-algebra $C(\times_\bn X)$ consisting of all continuous functions on the compact space $\times_\bn X$.

To pass immediately to the natural quantum generalisation, one can consider the infinite tensor product $\otimes_\bn\ga$
of an arbitrary (usually separable) unital $C^*$-algebra $\ga$.\footnote{Contrarily to the commutative case for which the $C^*$-cross norm in forming the tensor product is unique, in noncommutative case where such a cross norm is not uniquely determined in general, the minimal $C^*$-cross norm is usually used for such kind of questions, see {\it e.g.} \cite{St2}.}

On $\otimes_\bn\ga$ equipped with the norm arising from the minimal $C^*$-cross norm, the group of all finite permutations is also naturally acting, and the De Finetti theorem would characterise the structure of the set of invariant states under such an action when $\ga$ is abelian. This quantum extension was carried out by St{\o}rmer in \cite{St2}. Indeed, among many interesting results in the lines of operator algebras, in that paper it was proved that any symmetric state $\om$ is the "barycentre" ({\it i.e.} a mixture) of a measure $\n_\om$, supported on the so-called product states, the last being the natural quantum generalisation of a sequence of independent and identically distributed random variables, which are then extremal (or ergodic) w.r.t. the action of the permutations.

The investigation of dynamical systems based on $\bz^2$-graded $C^*$-algebras, hence describing the possible presence of Fermi particles, is a natural topic for the obvious applications to quantum field theory, see {\it e.g.} \cite{BR2}. Very recently, such an investigation was intensively extended to spin models on lattices whose site-variables enjoy the so-called Canonical Anti-commutation Relation (CAR for short). The early investigation of such discrete model was systematically started in \cite{AM1}, and continued in a long series of papers. We also mention the recent investigation in \cite{BF},
of disordered models ({\it i.e.} including the Fermi counterpart of a spin-glass) based on the $\bz^2$-graded $C^*$-algebras.

Recently in \cite{AM1, AM2}, the concept of the {\it product state} was extended to the so-called CAR algebra, which is the $C^*$-algebra describing Fermi particles ({\it i.e.} particles obeying to the Fermi statistics). It was also shown ({\it cf.} 
\cite{AFM, F1}) that the Fermi product state has an important role in investigating the Markov-like stochastic processes on the CAR algebra.

With the scope to extend and investigate the detailed balance for Fermi systems, the abstract ({\it i.e.} applicable to general $\bz^2$-graded $C^*$-algebras) concept of $C^*$-Fermi tensor product of two $\bz^2$-graded $C^*$-algebras was introduced and studied in \cite{CDF}. A crucial preliminary investigation in order to built such a Fermi tensor product was that involving the product states.

In the last mentioned paper, the completion of the algebraic part of the Fermi tensor product was carried out by using the maximal $C^*$-Fermi cross norm, which was crucial to define and investigate the so called diagonal, or maximally entangled, state. Yet, for the sake of completeness, also the the Fermi analogous of the minimal $C^*$-cross norm was introduced at the end of Section 8. Therefore, since there are now all the needed ingredient, the extension of the De Finetti theorem to general Fermi systems is a natural research line to be addressed.

The present note, which is the first part of a project devoted to the study of the structure and the main properties of the relative symmetric states, deals precisely with establishing the De Finetti theorem for general Fermi systems. It is a natural extension, first of the corresponding theorem for infinite $C^*$-tensor product of \cite{St2}, and then of that in \cite{CF} established for the particular case of the CAR algebra.

After some preliminary notions, we first describe in detail the Fermi $C^*$-tensor product $\ga\, \circled{\rm{{\tiny F}}}\, \gb$
between two $\bz^2$-graded $C^*$-algebras $\ga$ and $\gb$, built using the minimal norm firstly introduced in Section 8 of \cite{CDF}. Among other facts, we show that the natural grading on the algebraic part $\ga\, \circled{\rm{{\tiny F}}}_o\, \gb$ easily extends to a grading on the $C^*$-completion w.r.t. the minimal norm $\ga\, \circled{\rm{{\tiny F}}}\, \gb$, which makes the latter as a $\bz^2$-graded $C^*$-algebra in a natural way. 

The second step is that to construct the infinite Fermi $C^*$-tensor product $\ga=\circled{\rm{{\tiny F}}}_{\,\bn}\, \gb$ of a single $\bz^2$-graded $C^*$-algebra $\gb$, and notice that such an infinite product is still a $\bz^2$-graded $C^*$-algebra in a natural way. The group consisting of all finite permutation 
$\bp_\bn$ of $\bn$ acts in a natural way on such an infinite product $\ga$. 

After that, the paper continues with the investigation of the main ergodic useful properties of the dynamical system $(\ga, \bp_\bn)$. Among those, we show that the symmetric states ({\it i.e.} those invariant under the action of the permutations) is made of even ones w.r.t. the grading of $\ga$ induced by that of $\gb$, and 
$(\ga, \bp_\bn)$ is $\bp_\bn$-abelian. This allow us to conclude that
the symmetric states constitute a simplex made of even states.

The remaining step to establish the Fermi version of the De Finetti theorem is to show that the convex boundary is made of product states, each of them being a product of a single even state on $\gb$. This is proven in the last section of the present paper.

Some other interesting properties of the set of symmetric states are discussed in the forthcoming paper \cite{Fdue}.

\section{Preliminaries}

We gather here some facts useful for the forthcoming sections.
\medskip

\noindent
\textbf{Basic notions.} Let $\ch$ be a Hilbert space, and $\xi\in\ch$. The corresponding vector functional on $\cb(\ch)$ is defined as 
$\om_\xi:=\langle\,{\bf\cdot}\,,\xi,\xi\rangle$.

Let $\ga$ be a $C^*$-algebra. With $\aut(\ga)$ we denote the group of its $*$-automorphisms, with identity $e\aut(\ga)=\id_\ga$.

For two linear spaces $X$ and $Y$, we denote by $X\dot{+}Y$ and $X\odot Y$ their algebraic direct sum and tensor product, respectively.
If in addition, $\ga$ and $\gb$ are involutive algebras with $*$ denoting their involution, then
$\ga\otimes\gb$
will denote the algebraic tensor product $\ga\odot\gb$ equipped with the usual product and involution given on the generators by
$$
(a_1\otimes b_1){\bf\cdot}(a_2\otimes b_2)=a_1a_2\otimes b_1b_2\,,\quad (a_1\otimes b_1)^\dagger=a_1^*\otimes b_1^*\,,
$$
for all $a_1\,a_2\in\ga, b_1,b_2\in\gb\,.$\footnote{We introduce the symbols ${\bf\cdot}$ and ${}^\dagger$ to denote the product and the involution in
$\ga\otimes\gb$ in order to distinguish them from the analogous operations (denoted by the standard symbology)
in the $\bz_2$-graded tensor product $\ga\, \circled{\text{{\tiny F}}}\, \gb$ whenever $\ga$ and $\gb$ are equipped with a $\bz_2$-grading.}


For $C^*$-algebras $\ga$ and $\gb$, $\ga\otimes_{\rm max}\gb$
and $\ga\otimes_{\rm min}\gb$ are the completion of $\ga\otimes\gb$ w.r.t. the maximal and minimal $C^*$-cross norm, respectively, see {\it e.g.} \cite{T}.

For states $\om\in\cs(\ga)$, $\f\in\cs(\gb)$ on the $C^*$-algebras $\ga$ and $\gb$, we denote by
$\psi_{\om,\f}\in\cs(\ga\otimes_{\rm min}\gb)$
the product state on the $C^*$-algebra $\ga\otimes_{\rm min}\gb$.
A fortiori, $\psi_{\om,\f}$ is also well defined as a state on $\ga\otimes_{\rm max}\gb$. Therefore, with an abuse of notation we write $\psi_{\om,\f}\in\cs(\ga\otimes_{\rm min}\gb)$, or merely on its algebraic part
$\psi_{\om,\f}\in\cs(\ga\otimes\gb)$.

Let $\ga$ be a $C^*$-algebra, and $\f\in\cs(\ga)$ a state.
By $\big(\ch_\f,\pi_\f,\xi_\f\big)$, we denote the Gelfand-Naimark-Segal (GNS for short) representation associated to the state $\f$, see {\it e.g. } \cite{T}. If in addition $\th\in\aut(\ga)$ is a $*$-automorphism leaving invariant the state $\f$, then there exists a unitary $V_{\f,\th}$ acting on $\ch_\f$ which implements $\th$, that is
$$
V_{\f,\th}\pi_\f(a)V_{\f,\th}^*=\pi_\f(\th(a))\,,\quad a\in\ga\,.
$$

The quadruple $\big(\ch_\f,\pi_\f, V_{\f,\th},\xi_\f\big)$ is called {\it the covariant GNS representation} associated to the invariant state $\f$.
When the involved automorphisms $\{\a_g\mid g\in G\}\subset\aut(\ga)\}$ come from an action $G\ni g\mapsto\a_g\in \aut(\ga)$ and the state $\f$ is invariant under all $\a_g$, we denote such a group of unitaries in $\cb(\ch_\f)$ by $\{U_\f(g)\mid g\in G\}$.



\medskip

\noindent
\textbf{$\bz^2$-graded $C^*$-algebras.}





In the context of the present paper, a $\bz_2$-graded $C^*$-algebra is a pair made $(\ga,\th)$ of a $C^*$-algebra and a $*$-automorphism $\th$ of $\ga$ such that 
$\th^2=e_{\aut(\ga)}=\id_\ga$. Sometimes, we also consider $\bz_2$-grading on merely involutive algebras, see \cite{CDF}, Section 4.

It is immediate that any $a\in\ga$ can be written in a unique way as 
$$
a=\frac{a+\th(a)}2+\frac{a-\th(a)}2=:a_++a_-\,,
$$
where the corresponding {\it grade} is $\partial(a_{\pm})=\pm1$. Therefore, 
$$
\ga=\frac{(\id_\ga+\th)(\ga)}2\dot{+}\frac{(\id_\ga-\th)(\ga)}2=:\ga_++\ga_-\,,
$$
with $\ga_+\subset\ga$ is a $C^*$-subalgebra, and $\ga_-\subset\ga$ is a closed linear subspace. In addition
\begin{equation*}
\th\lceil_{\ga_1}=\id_{\ga_1}\,,\quad \th\lceil_{\ga_{-1}}=-\id_{\ga_{-1}}\,,
\end{equation*}
and the elements of $\ga_+$ and $\ga_-$ are the homogeneous parts of $\ga$, and are called the {\it even} and the {\it odd} part of $\ga$, respectively. 

The map $E:=(\id_\ga+\th)/2:\ga\to\ga_+$ is indeed a conditional expectation onto $\ga_+$ satisfying $E\circ\th=E=\th\circ E$.

Sometimes, we omit to indicate the grading automorphism $\th$ when this causes no confusion.

Let $\big(\ga^{(i)},\th^{(i)}\big)$, $i=1,2$, be a pair of $\bz_2$-graded $C^*$-algebras,
together with a map $T:\ga^{(1)}\to\ga^{(2)}$. $T$ is said to be {\it even} if it is grading-equivariant:
$$
T\circ\th^{(1)}=\th^{(2)}\circ T\,.
$$

When $\th^{(2)}=\id_{\ga^{(2)}}$, the map $T:\ga^{(1)}\to\ga^{(2)}$ is even if and only if it is grading-invariant, that is
$T\circ\th^{(1)}=T$. As a particular case when $\big(\ga^{(1)},\th^{(1)}\big)=(\ga,\th)$
and $\big(\ga^{(2)},\id_{\ga^{(2)}}\big)=\big(\bc,\id_\bc\big)$, a functional $f:\ga\to\bc$
is even if and only if $f\circ\th=f$. If the map $T$ (or the functional $f$) is $\bz_2$-linear,
then it is even if and only if $T\lceil_{\ga^{(1)}_-}=0$ ($f\lceil_{\ga_-}=0$).

Particular importance has the set of the {\it even states}, denoted by $\cs_+(\ga)$. By the previous considerations 
$$
\f\in\cs(\ga)\,\,\text{is even}\,\,\iff \f\lceil_{\ga_-}=0\,,
$$
and thus $\cs_+(\ga)\sim\cs(\ga_+)$ as topological spaces.

The most common example of nontrivial graded $C^*$-algebra, is the CAR algebra based on an infinite countable set, see {\it e.g.} \cite{T}, Exercise X1V.1.

\medskip

\noindent
\textbf{Dynamical systems.} A $C^*$-dynamical system is a triple $(\ga, G, \a)$ made of a $C^*$-algebra, a group $G$ and a representation, usually called an {\it action}, $G\ni g\mapsto\a_g\in\aut(\ga)$ of $G$ by $*$-automorphisms of $\ga$. 

If there is no matter of confusion, sometimes we omit to indicate the symbol $\a$, and denote such a dynamical system as $(\ga,G)$. For example, the fixed-point subalgebra, (which is a $C^*$-subalgebra of $\ga$) made of invariant element of $\ga$ under the action of $G$, is simply denoted by $\ga^G$. 

We denote by $\cs_G(\ga)\subset\cs(\ga)$ the locally compact, under the $*$-weak topology, convex set of the invariant states under the action of $G$. If $\ga$ is unital with identity $\idd_\ga=:\idd$, $\cs_G(\ga)$ is compact, and the convex boundary set $\partial\cs_G(\ga)$ made of the extremal invariant states will be denoted by $\ce_G(\ga)$ Such extremal invariant states are also called ergodic. From now on, we suppose that $\ga$ is unital, if it is not otherwise specified.

Let $\f\in\cs_G(\ga)$ with $(\pi_\f,\ch_\f,U_\f,\xi_\f)$ be the covariant GNS representation associated to $\f$. The action 
$\ad_{U_\f(g)}$ of $G$ on $\cb(\ch_\f)$ leaves invariant both $\pi_\f(\ga)''$ and $\pi_\f(\ga)'$, and hence its centre 
$\ty{Z}_\f:=\pi_\f(\ga)''\cap\pi_\f(\ga)'$. We denote the invariant elements under such an action, simply by $\ty{Z}_\f^G$. We have then
$$
{\rm Z}_\f^G=\pi_\f(\ga)''\cap\pi_\f(\ga)'\cap U_\f(G)'\,.
$$

For a state $\f\in\cs_G(\ga)$, we denote $E_\f$ the selfadjoint projection onto the (closed) subspace of $\ch_\f$ of the invariant vectors under
$U_\f(g)$, for all $g\in G$. The invariant state $\f$ is said to be $G$-{\it abelian} if the subspace $E_\f\pi_\f(\ga)E_\f$ is made of mutually commuting operators. The whole dynamical system $(\ga,G)$ is said to be $G$-abelian if any invariant state is.

For a state $\f\in\cs_G(\ga)$, consider the following assertions:
\begin{itemize}
\item[(i)] $\f$ is ergodic,
\item[(ii)] $\{\pi_\f(\ga),U_\f(G)\}'=\bc\idd_{\ch_\f}$,
\item[(iii)] $E_\f\ch_\f$ is one-dimensional.
\end{itemize}
Then ({\it e.g.} \cite{S}, Section 3.1) (iii)$\Longrightarrow$(ii)$\iff$(i), and those are all equivalent if $\f$ is $G$-abelian.

The state $\f$ is said to be {\it weakly clustering} if
$$
M\{\f(a\a_{\hat g}(b)\}=\f(a)\f(b)\,,\quad a,b\in\ga\,,
$$
where $M$ is any mean obtained by using a F\o lner sequence if $G$ is amenable.\footnote{Here, "$\,\hat { }\,\,$'' denotes a dumb variable.}

It is possible to see that (iii) above is equivalent to $\f$ being weakly clustering, see {\it e.g.} 
\cite{F22}.\footnote{The analysis can be extended to any group $G$ by using the Godement mean, see \cite{DK}.}

A state $\f\in\cs_G(\ga)$ is said to be {\it strongly
clustering} if there exists a sequence $(g_n)_n\subset G$ such that, for every $a,b\in\ga$,
$$
\lim_n\f(\a_{g_n}(a)b)=\f(a)\f(b)\,.
$$
 
It is expected that $\f$ being strongly clustering implies $\f$ being weakly clustering, but the reverse does not hold. We never use such an expected result, but it is shown for the infinite $C^*$-tensor product and for the infinite CAR algebra ({\it cf.} \cite{St2} and \cite{CF}, respectively)
that a weakly clustering symmetric state is strongly clustering. We will see that this holds true (essentially under the same proof as in \cite{CF}) also for the infinite $C^*$-Fermi tensor product of any (unital) $C^*$-algebra.

We say that $G$ is acting on $\ga$ as a {\it large group of automorphisms} if, for any $\f\in\cs_G(\ga)$,
\begin{align*}
&\overline{{\rm conv}(\pi_\f(\a_g(a):g\in G)}\cap\pi_\f(\ga)'\neq\emptyset\,,\\
&\qquad\quad\quad\text{for each selfadjoint}\, a\in\ga\,.
\end{align*}

It is shown in \cite{St1}, Theorem 3.1 that, if $\f\in\cs_G(\ga)$, then there exists a unique normal $G$-invariant conditional expectation
$\F_\f:\pi_\f(\ga)''\to\ty{Z}^G_\f$ of $\pi_\f(\ga)''$ onto $\ty{Z}^G_\f$.

\section{Tensor product of Fermi algebras}


For $\bz_2$-graded $C^*$-algebras $(\ga,\a)$, $(\gb,\b)$,
we denote by 
$\ga\, \circled{\rm{{\tiny F}}}_o\, \gb$ the algebraic tensor product 
$\ga\odot\gb$ equipped with the $*$-operation and product as follows.

One first notices that, as linear spaces, 
 $$
\ga\, \circled{\rm{{\tiny F}}}_o\, \gb=\ga\odot\gb=\ga\otimes\gb\,,
$$
the last being the algebraic tensor product equipped with its natural "star"-operation denoted by "$\dagger$" and the product denoted by a "dot". 
Therefore,
$$
\ga\odot\gb=\dot{+}_{i,j\in\bz_2}(\ga_i\odot\gb_j)=\ga\otimes\gb\,.
$$

For homogeneous elements $a\in\ga$, $b\in\gb$ and $i,j\in\bz_2$, we recall the following definition
\begin{align*}
\eeps(a,b):=\left\{\!\!\!\begin{array}{ll}
                      -1 &\text{if}\,\, \partial(a),\partial(b)=-1\,,\\
                     \,\,\,\,\,1 &\text{otherwise}\,.
                    \end{array}
                    \right.
\end{align*}

Consider the generic elements $x,y\in\ga\odot\gb$. We can write
\begin{align*}
\begin{split}
&x:=\oplus_{i,j\in\bz_2}x_{i,j}\in\oplus_{i,j\in\bz_2}(\ga_i\odot\gb_j)\,,\\
&y:=\oplus_{i,j\in\bz_2}y_{i,j}\in\oplus_{i,j\in\bz_2}(\ga_i\odot\gb_j)\,,
\end{split}
\end{align*}
and we set
\begin{equation}
\label{prstc}
\begin{split}
x^*:=&\sum_{i,j\in\bz_2}\eps(i,j)x_{i,j}^\dagger\, ,\\
xy:=&\sum_{i,j,k,l\in\bz_2}\eps(j,k)x_{i,j}{\bf\cdot}y_{k,l}\,.
\end{split}
\end{equation}
Notice that, for $x=a\odot b$ and $y=A\odot B$, where $a,A\in \ga$ and $b,B\in\gb$,
\begin{align*}
x^*=&\eeps(a,b)x^\dagger=\eeps(a,b)a^*\odot b^*\,,\\
xy=&\eeps(b,A)(a\otimes b){\bf\cdot}(A\otimes B)=\eeps(b,A)aA\otimes bB\,.
\end{align*}

Under the above operations ({\it cf.} \cite{CDF}, Proposition 6.1), the involution and product operations on $\ga\odot\gb$ given in \eqref{prstc} are well defined, and make the linear space $\ga\odot\gb$ an involutive algebra, indeed denoted by $\ga\, \circled{\rm{{\tiny F}}}_o\, \gb$.

As noted in Section 6 of  \cite{CDF}, for involutive $\bz_2$-graded ($C^*$-)algebras $(\ga,\a)$ and $(\gb,\b)$,
their Fermi tensor product $\ga\, \circled{\rm{{\tiny F}}}_o\, \gb$ is naturally equipped with a structure of involutive
$\bz_2$-graded algebra, by putting
\begin{equation}
\label{picm}
\begin{split}
\big(\ga\, \circled{\rm{{\tiny F}}}_o\, \gb\big)_+:=&\big(\ga_+\odot\gb_+\big)\dot{+}\big(\ga_{-}\odot\gb_{-}\big)\,,\\
\big(\ga\, \circled{\rm{{\tiny F}}}_o\, \gb\big)_-:=&\big(\ga_+\odot\gb_{-}\big)\dot{+}\big(\ga_{-}\odot\gb_{+}\big)\,.
\end{split}
\end{equation}

In this situation, such a grading is induced by the involutive automorphism $\th_o=\a\, \circled{\rm{{\tiny F}}}\, \b$ given on the elementary tensors by
\begin{equation}
\label{aupf}
(\a\, \circled{\rm{{\tiny F}}}\, \b)(a\, \circled{\rm{{\tiny F}}}\, b):=\a(a)\, \circled{\rm{{\tiny F}}}\, \b(b)\,,\quad a\in\ga\,,\,\, b\in\gb\,.
\end{equation}

We now report the main properties of the so-called {\it product state}. Indeed for the involutive $\bz_2$-graded algebras $(\ga,\a)$ and $(\gb,\b)$, and for given $\om\in\cs(\ga)$ and $\f\in\cs(\gb)$, the first step is to understand whether the product functional of $\om$ and $\f$ is well defined and
positive on the algebraic Fermi tensor product $\ga\, \circled{\rm{{\tiny F}}}_o\, \gb$.

If needed, we suppose without loosing generality that $\ga$ and $\gb$ are
unital with units $\idd_\ga$ and $\idd_\gb$, respectively. The general situation can be
achieved by adding the identities, or considering approximate identities on $\ga$ and $\gb$ as well.

We now report the main properties of the product functional, firstly studied in \cite{AM1, AM2} for concrete systems based
on the CAR algebra, and then in Section 7 of \cite{CDF} for abstract Fermi systems.

Indeed, for $\om\in\cs(\ga)$ and $\f\in\cs(\gb)$, the product functional $\om\times\f$ on $\ga\, \circled{\rm{{\tiny F}}}\,\gb$ is defined as usual by
\begin{equation}
\label{finprc}
\om\times\f\bigg(\sum_{j=1}^n a_j\, \circled{\rm{{\tiny F}}}\,b_j\bigg):=\sum_{j=1}^n\om(a_j)\f(b_j),\quad \sum_{j=1}^n a_j\, \circled{\rm{{\tiny F}}}\,b_j\in \ga\, \circled{\rm{{\tiny F}}}\,\gb\,.
\end{equation}

It is well known
that the functional given in \eqref{finprc} is well defined on $\ga\odot\gb$, and
therefore on $\ga\, \circled{\rm{{\tiny F}}}_o\, \gb$. In addition, even though it coincides
with the corresponding product state on $\ga\otimes \gb$, denoted by $\psi_{\om,\f}$, n general it is not positive on the involutive algebra 
$\ga\, \circled{\rm{{\tiny F}}}\, \gb$, and the necessary and sufficient condition for positivity is that at least one of $\om$ or $\f$ must be even. This is the content of Proposition 7.1 in \cite{CDF}, in which it was proven the crucial facts that, for each $x,y\in\ga\, \circled{\rm{{\tiny F}}}_o\, \gb$,
\begin{equation*}
\begin{split}
|(\om\times\f)(x)|\leq&(\om\times\f)(x^*x)^{1/2}\,,\\
(\om\times\f)(x^*y^*yx)\leq&C_y(\om\times\f)(x^*x)\,,
\end{split}
\end{equation*}
where $C_y$ is a positive constant depending on $y$.

By the above properties ({\it e.g.} \cite{CDF}, Theorem 3.2), the GNS representation of $\f\times\psi$ provides operators $\pi_{\f\times\psi}$ which are bounded, and the product state itself is given by a (unique, up to a multiplicative phase) vector state $\om_{\xi_{\om\times\f}}$ associated to the  cyclic vector $\xi_{\om\times\f}$.

For the purpose of the present note, we use the completion of $\ga\, \circled{\rm{{\tiny F}}}_o\, \gb$ under the norm $[]\,\,\,[]$ introduced at the end of Section 8 of \cite{CDF}. Indeed, define
\begin{equation*}
[]c[]_{\rm min}:=\sup\{\|\pi_{\f\times\psi}(c)\|\mid \f\in\cs(\ga)_+,\,\,\psi\in\cs(\gb)_+\},\quad c\in\ga\, \circled{\rm{{\tiny F}}}_o\,\gb\,.
\end{equation*}
In  \cite{CDF}, it was shown that $[]\,\,[]_{\rm min}$ is indeed a norm on $\ga\, \circled{\rm{{\tiny F}}}_o\,\gb$, which can be viewed as the analogous of the minimal $C^*$-cross norm for $\ga\otimes\gb$.\footnote{ In \cite{CRZ}, it was shown that such a norm is indeed minimal among all the $C^*$-cross norms on $\ga\, \circled{\rm{{\tiny F}}}_o\,\gb$.}

We also note that the maximal norm $[]\,\,[]_{\rm max}$, introduced and studied in \cite{CDF}, plays an important role in ergodic theory involving abstract $C^*$-Fermi systems. We mention the extension of the concept of diagonal state, being orthogonal to the corresponding product state in most of interesting cases, where the maximal norm plays a crucial role.

Since we use only the minimal norm, we put $\ga\, \circled{\rm{{\tiny F}}}\, \gb:=\overline{\ga\, \circled{\rm{{\tiny F}}}_o\,\gb}^{[]\,\,[]_{\rm min}}$, always denoting such a minimal norm simply by $\|\,\,\,\|$.
\begin{rem}
\label{unprid}
We point out that, the GNS representation $\pi_{\f\times\psi}$ and the product state $\f\times\psi$ uniquely extend to the whole
$\ga\, \circled{\rm{{\tiny F}}}\, \gb$ to a representation and a state, denoted again by $\pi_{\f\times\psi}$ and $\f\times\psi$ respectively.
\end{rem}
\begin{proof}
Indeed, for $c\in\ga\, \circled{\rm{{\tiny F}}}_o\,\gb$, we easily have $\|\pi_{\f\times\psi}(c)\|_{\cb(\ch_{\f\times\psi})}\leq\|c\|$, and thus $\pi_{\f\times\psi}$ uniquely extend to a contractive map of the completion $\ga\, \circled{\rm{{\tiny F}}}\,\gb$ of $\ga\, \circled{\rm{{\tiny F}}}_o\,\gb$ into $\cb(\ch_{\f\times\psi})$ which is indeed a 
representation. 

Concerning the product state, it is uniquely extended to the whole $\ga\, \circled{\rm{{\tiny F}}}\,\gb$ by the vector state $\om_{\xi_{\f\times\psi}}$.\end{proof}

If $\f\in\ga$, $\psi\in\gb$, first we remark that $\f\times\psi$ is always well defined as linear functional on $\ga\, \circled{\rm{{\tiny F}}}_o\,\gb$. By \cite{CDF}, Proposition 7.1, it is positive if and only if at least one of them is even. By the above consideration, its GNS representation and itself uniquely extend to the whole $\ga\, \circled{\rm{{\tiny F}}}\,\gb$. 

Therefore, we can conclude that the product state $\f\times\psi$, which is indeed a state on 
$\ga\, \circled{\rm{{\tiny F}}}\,\gb$, is uniquely determined by its values on the generators of $\ga\, \circled{\rm{{\tiny F}}}\,\gb$, provided it exists.

\section{The infinite Fermi tensor product}

We show that the completion under the minimal norm of the infinite tensor product is indeed a $\bz_2$-graded $C^*$-algebra, a fact already implicitly noted in \cite{CDF}.

Let $\ga\, \circled{\rm{{\tiny F}}}_o\, \gb$ be the algebraic Fermi product of $\bz_2$-graded $C^*$-algebras $(\ga,\a)$, $(\gb,\b)$, equipped with the grading $\th_o:=\a\, \circled{\rm{{\tiny F}}}_o\, \b$ defined by \eqref{picm} and \eqref{aupf}.
\begin{prop}
\label{rgex}
The involutive grading $*$-automorphism $\th_o$ uniquely extends to an involutive $*$-automorphism $\th$ of the whole $\ga\, \circled{\rm{{\tiny F}}}\, \gb$, and thus makes $(\ga\, \circled{\rm{{\tiny F}}}\, \gb,\th)$ as a $\bz_2$-graded $C^*$-algebra.
\end{prop} 
\begin{proof}
Fix $\f\in\cs_+(\ga)$ and $\psi\in\cs_+(\gb)$, and denote by $\om:=\f\times\psi$ the corresponding product state on $\ga\, \circled{\rm{{\tiny F}}}_o\, \gb$. Define
$$
V_o\pi_\om(a)\xi_\om:=\pi_\om(\th_o(a))\xi_\om\,,\quad a\in\ga\, \circled{\rm{{\tiny F}}}_o\, \gb\,.
$$

Since $\f$ and $\psi$ are even, $\om$ is $\th_o$-invariant, and thus
$$
\|V_o\pi_\om(x)\xi_\om\|^2=\om(\th_o(x^*x))=\om(x^*x)=\|\pi_\om(x)\xi_\om\|^2\,.
$$

Therefore $V_o$ is a well-defined isometry with dense range whose closure is a unitary $V$ acting on $\ch_\om$ leaving $\xi_\om$ invariant and satisfying 
$$
V\pi_\om(x)V^*=\pi_\om(\th_o(x))\,,\quad x\in\ga\, \circled{\rm{{\tiny F}}}_o\, \gb\,,
$$
which implies that $\|\pi_\om(\th_o(a))\|=\|\pi_\om(a)\|$.

We then conclude that $\th_o$ is isometric on $\ga\, \circled{\rm{{\tiny F}}}_o\, \gb$ w.r.t. the minimal norm, and thus it extends to the whole 
$\ga\, \circled{\rm{{\tiny F}}}_o\, \gb$ to an isometric map, denoted by $\a\, \circled{\rm{{\tiny F}}}\, \b$, which is indeed an involutive $*$-automorphism.
\end{proof}

According to the purpose of the present note, from now on, we consider a unital $\bz^2$-graded $C^*$-algebra $(\gb,\a)$ with the unique, necessarily even identity $\idd$, together with a discrete countable set which, in our situation, can be chosen as $\bn$, and form the $\bz^2$-graded (minimal) infinite $C^*$-Fermi tensor product 
$\big(\circled{\rm{{\tiny F}}}_{\,\bn}\, \gb, \circled{\rm{{\tiny F}}}_{\,\bn}\, \a\big)$ as follows.\footnote{We point out that the construction of the infinite Fermi tensor product can be straightforwardly carried out also for index-sets of arbitrary cardinality and non unital $C^*$-algebras.}

Indeed, we take into account Proposition \ref{rgex}, pointing out that, for each $n\in\bn$ there is a natural embedding 
$$
a\in\underbrace{\gb\, \circled{\rm{{\tiny F}}}\, \cdots \, \circled{\rm{{\tiny F}}}\,\gb}_{\textrm{n-times}}\mapsto
a\, \circled{\rm{{\tiny F}}}\, \idd
\in\underbrace{\gb\, \circled{\rm{{\tiny F}}}\, \cdots \, \circled{\rm{{\tiny F}}}\,\gb}_{\textrm{n-times}}\, \circled{\rm{{\tiny F}}}\,\gb
$$
between $\bz_2$-graded $C^*$-algebras,
which preserves the parity in the obvious way:
$$
\big(\underbrace{\a\, \circled{\rm{{\tiny F}}}\, \cdots \, \circled{\rm{{\tiny F}}}\,\a}_{\textrm{n-times}}\, \circled{\rm{{\tiny F}}}\,\a\big)
(a\, \circled{\rm{{\tiny F}}}\, \idd)
=\big(\underbrace{\a\, \circled{\rm{{\tiny F}}}\, \cdots \, \circled{\rm{{\tiny F}}}\,\a}_{\textrm{n-times}}\big)(a)\, \circled{\rm{{\tiny F}}}\,\idd\,.
$$

We then have the following
\begin{prop}
The direct limit  
$$
(\ga_o,\th_o):=\stackrel[\longrightarrow]{}{\lim}{}_{n\to\infty}
\big(\circled{\rm{{\tiny F}}}_{i=0}^n\, \gb, \circled{\rm{{\tiny F}}}_{i=0}^n\, \a\big)
$$
is an involutive $\bz^2$-graded algebra admitting a unique $C^*$-completion $(\ga,\th)$, which is the $C^*$-inductive limit of the sequence
$\big(\circled{\rm{{\tiny F}}}_{i=0}^n\, \gb, \circled{\rm{{\tiny F}}}_{i=0}^n\, \a\big)_n$ of compatible $\bz^2$-graded $C^*$-algebras.
\end{prop}
\begin{proof}
It is a direct consequence of the above considerations and Section L2 of \cite{WO}.
\end{proof}
In the present note, we call the above $\bz^2$-graded $C^*$-algebra $(\ga,\th)=:\big(\circled{\rm{{\tiny F}}}_{\,\bn}\, \gb, \circled{\rm{{\tiny F}}}_{\,\bn}\, \a\big)$ simply {\it the infinite Fermi tensor product} of the single $\bz^2$-graded $C^*$-algebra $(\gb,\a)$.

The elements of $\ga_o$, that is those for which are asymptotically tensorised by the identity, are called {\it localised elements}.

It is matter or routine to verify that, given a collection $\{\f_n\mid n\in\bn\}\subset\cs_+(\gb)$ of even states, the product state
$\om:=\times_{n\in\bn}\f_n$ is well defined as an even state on the infinite fermi tensor product.
A crucial importance for the purpose of the present note has the infinite product state $\times_{\bn}\f$ of a 
single one $\f\in\cs_=(\gb)$.

We end the present section by noticing that $\ty{CAR}(\bn)$ can be viewed as an infinite $C^*$-tensor product, again on $\bn$, of 
$\gb\sim\bm_{2^d}(\bc)$, $d\geq1$, by regrouping the local algebras $d$ at a time. In such a situation,
$$
\ty{CAR}(\bn)\sim\underbrace{\bm_2(\bc)\vee\cdots\vee\bm_2(\bc)}_{\textrm{d-times}}
\, \circled{\rm{{\tiny F}}}\,\underbrace{\bm_2(\bc)\vee\cdots\vee\bm_2(\bc)}_{\textrm{d-times}}
\, \circled{\rm{{\tiny F}}}\,\cdots\,,
$$
where, for this case for which the single algebra $\gb$ is finite dimensional, the norm on the $C^*$-Fermi product is uniquely determined.

Therefore, by regrouping the algebras, the De Finetti theorem could be established using the results in \cite{CF} for "single site" graded algebras 
$(\gb,\b)\sim(\ty{CAR}(2^d),\a)$, $d>1$, where $\a$ is the $\bz^2$-grading on $\ty{CAR}(2^d)$.

\section{Ergodic properties}

The main ingredient of the present note is the dynamical system $(\ga,\bp)$ made of all infinite Fermi tensor product $(\ga,\th)=\big(\circled{\rm{{\tiny F}}}_{\,\bn}\, \gb, \circled{\rm{{\tiny F}}}_{\,\bn}\, \a\big)$, on which the (discrete) group of the finite permutations $\bp:=\bp_\bn$ acts in a natural way. We also omit to indicate the underlying gradings when this cause no confusion.

Indeed, on the generators of the direct limit
$a=a_0\, \circled{\rm{{\tiny F}}}\, \cdots \, \circled{\rm{{\tiny F}}}\,a_k\, \circled{\rm{{\tiny F}}}\,\cdots$, such an action $\a^{(o)}_g:\ga_o\to\ga_o$ is defined as
$$
\a^{(o)}_g(a_0\, \circled{\rm{{\tiny F}}}\, \cdots \, \circled{\rm{{\tiny F}}}\,a_k\, \circled{\rm{{\tiny F}}}\,\cdots):=
\a_{g(0)}\, \circled{\rm{{\tiny F}}}\, \cdots \, \circled{\rm{{\tiny F}}}\,\a_{g(k)}\, \circled{\rm{{\tiny F}}}\,\cdots)\,,
$$
where $a$ is a localised element made of elementary tensors.

Such an action extends to an action, denoted by $\a$, of $\bp$ on the whole infinite Fermi tensor product $\ga$, which is grading-equivariant: 
$\a_g\circ \th=\th\circ\a_g$, for each $g\in\bp$.

In order to investigate the ergodic properties of symmetric states ({\it i.e.} those invariant under the action of the permutation group $\bp$), we start with the following
\begin{thm}
\label{mainer}
For each state $\om\in\cs_\bp(\ga)$,
\begin{itemize}
\item[(i)] $\om$ is even;
\item[(ii)] $\om$ is $\bp$-abelian;
\item[(iii)] $\om$ is asymptotically abelian in average:
$$
M\{\om(c[\a_{\hat g}(a),b]d)\}=0\,,\quad a,b,c,d\in\ga\,.
$$
\end{itemize}
\end{thm}
\begin{proof}
We start by fixing unit vectors $\xi,\eta\in\ch_\om$.

(i) By a standard approximation argument, we can reduce the matter to localised odd elements and show that $\om\in\cs_\bp(\ga)$ vanishes to each of such elements. Fix a norm one localised odd element $a\in\ga_-$, with $k>0$ some (finite) integer such that 
$a=c\, \circled{\rm{{\tiny F}}}\,\idd \, \circled{\rm{{\tiny F}}}\,\idd\, \circled{\rm{{\tiny F}}}\,\cdots$, with 
$c\in\underbrace{\gb\, \circled{\rm{{\tiny F}}}\, \cdots \, \circled{\rm{{\tiny F}}}\,\gb}_{\textrm{k-times}}$. Denote also
${\bf n}:=\{0,\dots,n\}$ the finite set made of exactly the first $n+1$
elements. 

By taking into account the $\a_g(a)$ and $\a_g(a^*)$ are localised and odd, and using first the von Neumann ergodic theorem \cite{CF}, Proposition 3.1, and then the estimate in \cite{CF}, 
Lemma 3.3, we have
\begin{align*}
|\langle\{E_\om\pi_\om(a)E_\om,E_\om\pi_\om(a^*)&E_\om\}\xi,\eta\rangle|
=\big|M\{\langle E_\om\pi_\om(\{a,\a_{\hat g}(a^*)\})E_\om\xi,\eta\rangle\}\big|\\
\leq&\lim_{n\to\infty}\bigg(\frac1{n+1}\frac{\left|\{g\in\bp_ {{\bf n+1}}\mid\,{\bf k}\cap g{\bf
k}\neq\emptyset\}\right|}{n!}\bigg)\\
\leq&c(k)\lim_{n\to\infty}\frac1{n+1}=0\,.
\end{align*}
This implies that $E_\om\pi_\om(a)E_\om=0$ for $\a\in\ga_-$, and thus
$$
\om(a)=\om_{\xi_\om}(E_\om\pi_\om(a)E_\om)=0\,.
$$

Since $\overline{\th_o(\ga_o)}=\th(\ga)$, and thus $\overline{(\ga_o)_-}=\ga_-$, for each odd $a$ with $\|a\|=1$ and $0<\eps<1$, there exists an odd $a_\eps$ with $\|a_\eps\|=1$ such that $\|a-a_\eps\|<\eps/4$, and thus $\|a^*-a^*_\eps\|<\eps/6$. An easy calculation leads to
$$
\big|\langle(\{E_\om\pi_\om(a)E_\om,E_\om\pi_\om(a^*)E_\om\}-\{E_\om\pi_\om(a_\eps)E_\om,E_\om\pi_\om(a_\eps^*)E_\om\})\xi,\eta\rangle\big|
<\eps\,,
$$
and thus the assertion follows as $\eps$ is arbitrarily close to zero.

(ii) By approximating elements with norm one localised ones $a,b\in\ga_o$ as before, and taking into account that $E_\om\pi_\om(c)E_\om$ vanishes if $c$ is odd, we can suppose that $a,b$ are even and localised. Then, reasoning as in (i), we get
$$
\langle[E_\om\pi_\om(a)E_\om,E_\om\pi_\om(b)E_\om]\xi,\eta\rangle
=M\{\langle E_\om\pi_\om([a,\a_{\hat g}(b)])E_\om\xi,\eta\rangle\}=0
$$
since, if $a=c\, \circled{\rm{{\tiny F}}}\,\idd \, \circled{\rm{{\tiny F}}}\,\idd\, \circled{\rm{{\tiny F}}}\,\cdots$ and 
$b=d\, \circled{\rm{{\tiny F}}}\,\idd \, \circled{\rm{{\tiny F}}}\,\idd\, \circled{\rm{{\tiny F}}}\,\cdots$ with  
$c,d\in\underbrace{\gb\, \circled{\rm{{\tiny F}}}\, \cdots \, \circled{\rm{{\tiny F}}}\,\gb}_{\textrm{k-times}}$, we have as before
$$
\big|M\{\langle E_\om\pi_\om([a,\a_{\hat g}(b)])E_\om\xi,\eta\rangle\}\big|\leq c(k)\lim_{n\to\infty}\frac1{n+1}=0\,.
$$
The statement follows by a standard approximation argument as before.

Concerning (iii), and taking into account the previous computations, the result can be reached by splitting $a,c,d$ in their even and odd parts. 

Indeed, suppose first that $a$
is even. We get
\begin{align*}
&M\{\om(c\a_{\hat g}(a)bd)\}=M\{\om(\a_{\hat g}(a)cbd)\}\\
=&\langle E_\om\pi_\om(a)E_\om\pi_\om(cbd)E_\om\xi_\om,\xi_\om\rangle
=\langle E_\om\pi_\om(cbd)E_\om\pi_\om(a)E_\om\xi_\om,\xi_\om\rangle\\
=&M\{\om(cbd\a_{\hat g}(a))\}=M\{\om(cb\a_{\hat g}(a)d)\}\,.
\end{align*}

Suppose now $a$ is odd, and split $c$ and $d$ in its odd and even parts. We now get
\begin{align*}
M\{\om(c&\a_{\hat g}(a)bd)\}=\pm M\{\om(\a_{\hat g}(a)cbd)\}\\
=&\pm\langle E_\om\pi_\om(a)E_\om\pi_\om(cbd)\xi_\om,\xi_\om\rangle
=0\,.
\end{align*}

On the other hand, 
\begin{align*}
M\{\om(c&b\a_{\hat g}(a)d)\}=\pm M\{\om(cbd\a_{\hat g}(a))\}\\
=&\pm\langle\pi_\om(cbd)E_\om\pi_\om(a)E_\om\xi_\om,\xi_\om\rangle
=0\,.
\end{align*}

Collecting together, we obtain the assertion.
\end{proof}
The previous theorem immediately leads to some important consequences listed below.
\begin{cor}
A state $\om\in\ce_\bp(\ga)$ is ergodic if and only if it is weakly clustering.
\end{cor}
\begin{proof}
It follows combining (ii) of Theorem \ref{mainer} with Proposition 3.1.12 of \cite{S}, by taking into account that 
$\om$ being strongly clustering is equivalent to $\dim(E_\om\ch_\om)=1$.
\end{proof}
\begin{cor}
\label{sysa}
The set of the invariant states $\cs_\bp(\ga)$ is a simplex in the sense of Choquet (cf \cite{CM}).
\end{cor}
\begin{proof}
It immediately follows from (ii) of the above theorem, by taking into account of Theorem 3.1.14 in \cite{S}.
\end{proof}
The following result is a direct consequence of the last corollary, by taking into account \cite{S}, Theorem 3.1.14.
\begin{thm}
\label{sysa1}
Let $\om\in\cs_\bp(\ga)$. Then there exists a unique probability Radon measure $\n_\om$ on $\cs_\bp(\ga)$ such that
\begin{itemize}
\item[(i)] $\D\subset\cs_\bp(\ga)$ Baire set with $\D\cap\ce_\bp(\ga)=\emptyset$, then $\n_\om(\D)=0$;
\item[(ii)] $\om=\int_{\cs_\bp(\ga)}\f\di\n_\om(\f)$.
\end{itemize}
\end{thm}
The measure $\n_\om$ is maximal, in the sense explained in \cite{S}, pag. 123, among the measures satisfying (ii) above and is nothing else than the measure associated to 
$\{\pi_\om(\ga),U_\om(\bp)\}'\subset\pi_\om(\ga)'$, see \cite{S}, Proposition 3.1.5.

We also note that, since $\n_\om$ satisfies (i) above, it is said that such a measure is {\it pseudo-supported} on $\ce_\bp(\ga)$. When $\ga$ is separable, $\n_\om$ is indeed supported on $\ce_\bp(\ga)$ because any Baire set is also a Borel set, and thus
\begin{equation}
\label{bary}
\om=\int_{\ce_\bp(\ga)}\f\di\n_\om(\f)\,,
\end{equation}
where \eqref{bary} is the so-called {\it barycentric} decomposition of $\om\in\cs_\bp(\ga)$.

In the present paper, the forthcoming results are never used for the proof of De Finetti theorem. We are providing these for the sake of completeness, together with Proposition \ref{scaz} below.
 \begin{prop}
 \label{avlgu}
The group $\bp$ acts on $\ga$ as a large group of automorphisms.
\end{prop}
\begin{proof}
Fix an arbitrary integer $n\in\bn$,
$a,b_1,\dots,b_n,c\in\ga$ with $a,b_1,\dots,b_n$ selfadjoint, $\om\in\cs_{\bp}(\ga)$, and denote by $\om_c:=\om(c^*(\,{\bf\cdot}\,)c)$.

Now, by (iii) of Theorem \ref{mainer}, for each $\eps>0$ there
exists a finite set $F\in\bn$ such that, for any $k=1,\dots,n$,
$$
\bigg|\om\bigg(c^*\Big[\bigg(\sum_{g\in \bp_F}\a_g(a)/|\bp_F|\bigg),b_k\bigg]c\bigg)\bigg|
=\bigg|\frac1{|\bp_F|}\sum_{g\in
\bp_F}\om(c^*[\a_g(a),b_k]c)\bigg|<\eps\,.
$$

Therefore,
\begin{align*}
0\leq&\inf\big\{|\om_c([\tilde a,b_k])|\mid \tilde a\in{\rm conv}(\a_g(a); g\in\bp)\big\}\\
\leq&\inf\bigg\{\bigg|\om\bigg(c^*\bigg[\bigg(\sum_{g\in \bp_F}\a_g(a)/|\bp_F|\bigg),b_k\bigg]c\bigg)\bigg|
\,; F\,\,\text{finite subset of}\,\,\bn\bigg\}\\
=&0\,,\quad k=1,2,\dots n\,.
\end{align*}

The assertion now follows by \cite{St1}, Theorem 3.5.
\end{proof}
\begin{prop}
Let $\om\in\cs_{\bp}(\ga)$ and $a\in\ga$. Then
\begin{equation*}
\mathop{\rm w}-\lim_{F\uparrow\bn}
\frac1{|\bp_F|}\sum_{g\in
\bp_I}\pi_\om(\a_g(a))=\F_\om(\pi_\om(a))\,,
\end{equation*}
where $\F_\om:\pi_\om(\ga)''\to\ty{Z}^\bp_\om$ is the normal $\bp$-invariant conditional expectation
 of $\pi_\om(\ga)''$ onto $\ty{Z}^\bp_\om$ which exists by Proposition \ref{avlgu}, and the limit is meant in the weak operator topology.
\end{prop}
\begin{proof}
For each fixed $a\in\ga$, consider any cluster point of the net 
$\big(\frac1{|\bp_{F}|}\sum_{g\in
\bp_{F}}\pi_\om(\a_g(a))\big)_{F\subset\bn}$, 
where $F$ denotes
any finite subsets of $\bn$, which exists by compactness. 
It is matter of routine ({\it e.g.} \cite{CF}, Proposition 4.3) to check that any such cluster point is $\ad_{U_\om(\bp)}$-invariant.

By using the asymptotic Abelianess (iii) of Theorem \ref{mainer}, and arguing as in the proof
of Lemma 5.3 of \cite{St1}, one shows that the limit above is indeed in $\ty{Z}^\bp_\om$.
But $\ty{Z}^\bp_\om$ can contain at most one of such limit points ({\it cf.} \cite{St1}, Theorem 3.1)
since $\bp$ acts on $\ga$ as a large group of automorphisms as we have just seen. As a consequence, the
limit is precisely $\F_\om(\pi_\om(a))$.
\end{proof}

\section{De Finetti theorem}

The present section is concerned with
the characterisation of the extremal symmetric states, which provides the extension of De Finetti theorem to the Fermi models under consideration. 

We also point out that some results of \cite{CF} can be straightforwardly adapted to the general situation of the present note. To simplify, we suppose that the one-site algebra $\gb$ is separable, even if most of the results hold true even in non separable case, as well as for infinite index sets of arbitrary cardinality.

We preliminary report the definition of the
sequence of permutations $(g_n)_n\subset\bp$ in \cite{St2} given
by
\begin{equation}
\label{mixi}
g_{n}\left( k\right) :=\left\{
\begin{array}{ll}
2^{n}+k & \text{if }0\leq k\leq 2^{n}\,, \\
k-2^{n} & \text{if }2^{n}<k\leq 2^{n+1}\,,\\
k & \text{if }2^{n+1}<k\,,
\end{array}
\right.
\end{equation}
which plays a crucial role as in \cite{St2, CF}.

\begin{lem}
\label{scaz1}
If $\om\in\ce_{\bp}(\ga)$, then for each $a\in\ga$,
$$
\mathop{\rm w}-\lim_{n}\big(\pi_\om(\a_{g_n}(a))\xi_\om\big)=\om(a)\xi_\om\,.
$$
\end{lem}
\begin{proof}
By a standard approximation argument, we can reduce the matter to localised elements
$a\in\ga_o$.  Let $g\in\bp$. Then there exists $n_{a,g}$ such that
$n>n_{a,g}$ implies $\a_{gg_n}(a)=\a_{g_n}(a)$. 

Hence, any weak limit point of the sequence $\big(\pi_\om(\a_{g_n}(a))\xi_\om\big)_n\subset\ch_\om$,
which exists by compactness, is an invariant vector
under the action of $\bp$, that is it belongs to $E_\om\ch_\om$.

Let $\xi\in\ch_\om$ be one such cluster points, and $(n_k)_k\subset\bn$ be a subsequence of natural numbers such that
\begin{align*}
\xi=&{\rm w}-\lim_{k}U_\om(g_{n_k})\pi_\om(a)U_\om(g_{n_k})^{-1}\xi_\om\\
=&{\rm
w}-\lim_{k}\pi_\om(\a_{g_{n_k}}(a))\xi_\om\,.
\end{align*}

Since $\xi\in E_\om\ch_\om$, and $\om$ is ergodic and thus $\dim(E_\om\ch_\om)=1$, there exists a constant $c\in\bc$ such that $\xi=c\xi_\om$.
Consequently,
$$
c=\lim_{k}\langle\pi_\om(\a_{g_{n_k}}(a))\xi_\om,\xi_\om\rangle=\om(a)\,,
$$
and thus there exists only one cluster point $\om(a)\xi_\om$.
\end{proof}
For the sake of completeness, we prove the following result, analogous to the corresponding ones in \cite{St1, CF}.
\begin{prop}
\label{scaz}
Let $\om\in\ce_\bp(\ga)$, and $a\in\ga$. For the following limit in the weak operator topology, we have
\begin{equation}
\label{aavvrr}
\mathop{\rm w}-\lim_{n}
\pi_\om(\a_{g_n}(a))=\om(a)\idd_{\ch_\om}=\F_\om(\pi_\om(a))\,.
\end{equation}
\end{prop}
\begin{proof}
By a standard approximation argument, it is enough to take vectors $\xi,\eta\in\ch_\om$ of the form
$\xi=\pi_\om(b)\xi_\om$, $\eta=\pi_\om(c^*)\xi_\om$, and reduce the
matter to $a,b,c\in\ga_o$. Let $a=a_++a_-$, $b=b_++b_-$ the
split of $a$, $b$ into their even and odd parts. By Theorem \ref{defin} below, it is enough to assume that $\om$ is
strongly clustering. 

Since $\om$ is even, by using the
standard (anti)commutation relations, we get
\begin{align*}
\lim_n\langle\pi_\om&(\a_{g_n}(a))\xi,\eta\rangle=\lim_n\om(c\a_{g_n}(a)b)\\
=&\lim_n\om(cb\a_{g_n}(a_+))+\lim_n\om(cb_+\a_{g_n}(a_-))-\lim_n\om(cb_-\a_{g_n}(a_-))\\
=&\om(cb)\om(a_+)+\om(cb_+)\om(a_-)-\om(cb_-)\om(a_-)\\
=&\om(cb)\om(a_+)+\om(cb_+)\om(a_-)\\
=&\om(cb)\om(a_+)+\om(cb_+)\om(a_-)+\om(cb_-)\om(a_-)\\
=&\om(cb)\om(a)=\langle(\om(a)\idd_{\ch_\om})\xi,\eta\rangle\,.
\end{align*}
Therefore, the first equality in \eqref{aavvrr} is satisfied for each $a\in\ga$.

Since $\om$ is ergodic, $\ty{Z}^\bp_\om\subset\{\pi_\om(\ga),U_\om(G)\}'=\bc\idd_{\ch_\om}$ and, since $\bp$ acts as a large group of automorphisms,
$\mathop{\rm w}-\lim_{n}\pi_\om(\a_{g_n}(a))$ must coincide with $\F_\om(\pi_\om(a))$, provided the former exists.
\end{proof}

The following result is the key-point to prove the generalisation of the De Finetti theorem for infinite $C^*$-Fermi tensor product of $\bz^2$-graded $C^*$-algebras. It generalises Theorem 5.3 in \cite{CF}, and can be proved following the lines of the corresponding Theorem 2.7 in \cite{St2}.
\begin{thm}
\label{defin}
Let $\om\in\cs_{\bp}(\ga)$. Then the
following are equivalent:
\begin{itemize}
\item[(i)] $\om$ is ergodic,
\item[(ii)] $\om$ is strongly clustering,
\item[(iii)] ${\displaystyle\om=\prod_\bn\f}$ for some even state $\f\in\gb$.
\end{itemize}
\end{thm}
\begin{proof}
$(i)\Rightarrow(ii)$ Suppose $\f$ is ergodic, and take
$a,b\in\ga$. Then by Lemma \ref{scaz1}, we get
$$
\lim_n\om(a\a_{g_n}(b))=\lim_n\langle\pi_\om(\a_{g_n}(a))\xi_\om,\pi_\om(a^*)\xi_\om\rangle=\om(a)\om(b)\,,
$$
that is $\om$ is strongly clustering.

$(ii)\Rightarrow(i)$  Choose a vector $\xi\perp\xi_\om$ belonging
to $E_\om\ch_\om$, and fix $\eps>0$. 

Since $\xi_\om$ is cyclic, there exists
$b\in\ga$ such that $\|\xi-\pi_\om(b)\xi_\om\|<\eps/2$, and thus 
$$
|\om(b)|=|\langle\pi_\om(b)\xi_\om,\xi_\om\rangle|=|\langle(\pi_\om(b)\xi_\om-\xi),\xi_\om\rangle|<\eps/2
$$
because $\xi$ is orthogonal to $\xi_\om$.

Let now
$a\in\ga$ such that $\|\p_\om(a)\|\leq1$. Recalling that $\xi$ and $\xi_\om$ are both $U_\om(\bp)$-invariant and considering the sequence in \eqref{mixi}, we get
\begin{align*}
|\langle\xi,\pi_\om(a)\xi_{\om}\rangle|=&|\langle U_{\om}(g_n)\pi_\om(a^*)U_{\om}(g_n)^{-1}\xi,\xi_{\om}\rangle|\\
\leq&|\langle U_{\om}(g_n)\pi_\om(a^*)U_{\om}(g_n)^{-1}\pi_\om(b)\xi_\om,\xi_\om\rangle|+\eps/2\\
=&|\om(\a_{g_n}(a^*)b)|+\eps/2\,.
\end{align*}

Suppose now $\om$ is strongly clustering and, taking the limit for
$n\rightarrow\infty$ on both sides, we get
$$
|\langle\xi,\pi_\om(a)\xi\om\rangle|\leq|\om(a^*)||\om(b)|+\eps/2<\eps/2+\eps/2=
<\eps\,.
$$

Since $\eps>0$ is arbitrary and $\xi_\om$ is cyclic for
$\pi_\om(\ga)$, we get $\xi=0$, and thus $E_\om\ch_\om$ is one dimensional, which implies that (and it is
indeed equivalent, since $(\ga,\bp)$ is $\bp$-abelian, to) the ergodicity of $\om$.

$(iii)\Rightarrow(ii)$ is obvious.

$(ii)\Rightarrow(iii)$ For $j\in\bn$, we denote the embedding 
$$
\gb\ni a\mapsto\iota_j(a):=
\underbrace{\idd\, \circled{\rm{{\tiny F}}}\, \cdots \, \circled{\rm{{\tiny F}}}\,\idd}_{\textrm{(j+1)-times}}\, 
\circled{\rm{{\tiny F}}}\, a\, \circled{\rm{{\tiny F}}}\,\idd \cdots\in\ga\,.
$$

It is well known that the product state is uniquely determined by the product of the values of the state on the
generators, see Remark \ref{unprid}. Therefore, the proof proceeds following the same lines of Theorem 2.7 of \cite{St2}, and then of Theorem 5.3 of \cite{CF}. We report the details for the convenience of the reader.

Indeed, for each $n\in\bn$ and
$a_0,\dots,a_n\in\gb$, we must show that
\begin{equation}
\label{prod1bis}
\om(\io_0(a_0)\cdots \io_n(a_n))=\prod_{j=0}^n\f(a_j)\,,
\end{equation}
for some state $\f\in\cs_+(\gb)$. The proof now proceeds as in Theorem 2.7 of \cite{St2}. 

For $j\in\bn$, define $\f_j:=\f\circ\io_j$. Since $\f$ is symmetric, hence even ({\it cf.} Theorem \ref{mainer}), the $\f_j$ are even and coincide each other: $\f_i=\f_j=:\f$ for $i,j\in\bn$. Now
\eqref{prod1bis} can be achieved by an induction procedure. 

Indeed,
for $n=0$ it follows immediately, so we suppose it holds true till
$n$. Fix $\eps>0$ and choose $m>n$ so large such that, using first the induction hypothesis and then the strong clustering property,  
\begin{align}
\label{cinque}
&\bigg|\om(\io_0(a_0)\cdots \io_n(a_{n})\a_{g_m}(\io_{(n+1)}(a_{(n+1)})))-\prod_{j=0}^{n+1}\f(a_j)\bigg|\\
=|\om(\io_0(a_0)&\cdots\io_n(a_{n})\a_{g_m}(\io_{n+1}(a_{n+1})))-\om(\io_0(a_0)\cdots
\io_n(a_{n}))\f(a_n)| <\eps\nn
\end{align}

Choose now a permutation $g\in\bp_\bn$ such that
$g(j)=j$ if $0\leq j\leq n$ and $g(n+1) =g_{m}(n+1) $. Then
\begin{align*}
\om(\io_0(a_0)\cdots \io_{n+1}(a_{n+1}))=&\om(\a_{g}(\io_0(a_0)\cdots \io_{n+1}(a_{n+1})))\\
=&\om(\io_0(a_0)\cdots \io_{n}(a_{n})\a_{g_m}(\io_{n+1}(a_{n+1})))
\end{align*}
which, combined with \eqref{cinque}, leads to 
$$
\bigg|\om(\io_0(a_0)\cdots \io_{n+1}(a_{n+1}))-\prod_{j=0}^{n+1}\f(a_j)\bigg|<\eps\,.
$$

The assertion follows as $\eps>0$ is arbitrary.
\end{proof}

We are in position to establish the version of De Finetti Theorem for Fermi systems, obtaining that any symmetric state is the mixture \eqref{bary} of product states, being each of them the product of a single even state.
\begin{thm}[De Finetti Theorem for infinite Fermi $C^*$-tensor products]
Let $\ga:=\circled{\rm{{\tiny F}}}_{\,\bn}\, \gb$ be the infinite Fermi $C^*$-tensor product of a single separable algebra $\gb$, together with the $C^*$-dynamical system $(\ga,\bp)$ obtained by considering the natural action of all finite permutations $\bp$ on $\ga$.
\begin{itemize}
\item[(i)] The set of the symmetric states $\cs_\bp(\ga)$ is a Choquet simplex;
\item[(ii)] the set of the ergodic states $\om\in\ce_\bp(\ga)$ is made of product states by a single even state $\f\in\cs_+(\gb)$: $\om=\times_\bn\f$;
\item[(iii)] for each $\om\in\cs_\bp(\ga)$, there exists a probability Radon measure $\n_\om$ supported on the ergodic states $\ce_\bp(\ga)$, such that
$\om$ is the barycentre of $\n_\om$: $\om=\int_{\ce_\bp(\ga)}\psi\di\n_\om(\psi)$.
\end{itemize}
\end{thm}
\begin{proof}
(i) is nothing else than Corollary \ref{sysa}, whereas (iii) immediately follows by Theorem \ref{sysa1} because, being $\gb$ and thus $\ga$, separable, the $\S$-algebra generated by the Baire sets of $\cs_\bp(\ga)$ coincides with that generated by the Borel sets. Finally, (ii) is nothing but the previous Theorem \ref{defin}.
\end{proof}

We end by noticing that the property for the set of extremal states to be $*$-weakly
closed implies a nice result. Indeed, in \cite{B} it is shown
that a simplex with closed boundary is affinely isomorphic to the
probability Radon measures on a compact Hausdorff space. This is the content of Theorem 2.8 of \cite{St2} (same proof), which in our situation assumes the form
\begin{prop}
\label{boundary} 
The Choquet simplex $\cs_{\bp}(\ga)$
has a $*$-weakly closed boundary and is affinely isomorphic to
the probability Radom measures on the convex $*$-weakly compact set $\cs(\gb_+)$.
\end{prop}

We also note that, when $\gb$ is $\bm_2(\bc)$, generated by the odd matrix $\begin{pmatrix} 
	 0 &1 \\
	0 & 0\\
     \end{pmatrix}$, then $\bm_2(\bc)_+\sim\bc^2$, and thus $\cs(\bm_2(\bc)_+)\sim[0,1]$ as was seen in Lemma 2.2 of \cite{CF}. Hence,  
Proposition 5.6 in \cite{CF} is a particular case of Proposition \ref{boundary}.

\end{document}